# A new version of distributional chaos and the relations between distributional chaos in a sequence and other concepts of chaos

Hongbo Zeng

Department of Mathematical Sciences, Tsinghua University, Beijing 100084, China

**Abstract**

In this paper we consider relations between distributional chaos in a sequence with distributional chaos, (ω-chaos, R-T chaos, DC3, respectively). We give a sufficient condition and prove that the distributional chaos is equivalent to the distributional chaos in a sequence under this condition. Besides, we get that distributional chaos in a sequence and ω-chaos(R-T chaos, DC3, respectively)don't imply each other. Finally, we give a new definition of chaos, named DC2', which is similar to DC2, and show that for any integer N > 0, f is DC2' if and only if $f^N$ is also DC2'.

**Keywords**：distributional chaos in a sequence; ω-chaos; R-T chaos; DC2'

## 1. Introduction

The existence of chaotic behavior in deterministic systems has attracted researchers for many years. Since Li and Yorke first gave the definition of chaos with strict mathematical terms in 1975 (see [1]), the research on chaos has had a great influence on modern science including natural science and many humanities, such as economics, sociology and philosophy. We can say that the influence has covered almost all disciplines. In almost all fields relating to dynamical progress, there exists a chaotic phenomenon. The theory of chaos convinces scientists that a simple definite system can produce complicated properties and a complex system possibly follow a simple law. As the ultimate aim of the scientists is to clarify the essence of the complexity, the chaotic systems with irregularly complex dynamical behaviors naturally become one common subject. However, depending on different perspective and understanding(people from different fields try to describe the behaviors by providing definitions of chaos according to their understanding of the subject.), the various concepts of chaos have been given, such as Li–Yorke chaos, Schweizer–Smital chaos (also called distributional chaos; see [2]), Devaney chaos [3], ω-chaos[4], R-T chaos[5], etc. Among them, distributional chaos has some actual significance. Later, three mutually nonequivalent versions of distributional chaos of type 1-3(DC1–DC3) were considered[6]. And more and more researchers give their attention to the properties of distributional chaos. Each definition tries to describe some kind of unpredictability of the system. Therefore, there exists much ambiguity in academic intercourse of different fields. At the same time, this situation cannot be tolerated in mathematical field which is based on strict mathematical definitions. Therefore, it is very significant to further explore the essence of chaos, unify the definition of chaos, and discuss the inner relations between the different definitions of chaos. Hence, in order to establish a satisfactory definitional and terminological framework for

E-mail address: hongbozeng92@163.com

chaotic system, we need to reveal the inner link between the various concepts which characterize the complexity. Since then, surely, it is an important question to understand the relation among the various definitions. And there are many results about that [6, 8-19]. In order to reveal the links between distributional chaos and chaos in the sense of Li and Yorke, the concept of distributional chaos in a sequence was introduced in [7]. In particular, The more related results about the distributional chaos in a sequence is referred in[7-11,17].

In this paper we study some different definitions of chaos and relations among them. We have known that the distributional chaos implies the distributional chaos in a sequence but the verse is not sure[11], in a word, the distributional chaos is stronger than the distributional chaos in a sequence. However, in this paper, we will give a sufficient condition to show that the distributional chaos is equivalent to the distributional chaos in a sequence under this condition. Besides, we give a new definition of chaos named distributional chaos of type 2'(DC2'),which contrast to DC2, and now the four versions of distributional chaos are mutually nonequivalent. Then we get a theorem that f is DC2' if and only if $f^N$ is too. Finally, we research the relationship about the other notions of chaos chaos( $\omega$-chaos, R_T chaos, DC3 ) with distributional chaos in a sequence.

This paper is organized as follows. In Section 2, we will first give some preliminaries and definitions. The main conclusions will be given in Section 3.

## 2. Preliminaries and definitions

We denote the $N$-fold iterates of f by $f^N$.

**Definition 1.** Let $(X, d)$ be a metric space. A continuous map $f: X \to X$ is called Li-Yorke chaotic if there exists an uncountable subset $S \subseteq X$ such that for every pair $x, y \in S$ of distinct points we have
$$\liminf_{n\to\infty} d(f^n(x), f^n(y)) = 0 \ and \ \limsup_{n\to\infty} d(f^n(x), f^n(y)) > 0.$$

**Definition 2.** Let $(X, d)$ be a metric space, and $f: X \to X$ be a continuous map. Note
$$F(f, x, y, t) = \liminf_{n\to\infty} \frac{1}{n} \#\{i \lceil d(f^i(x), f^i(y)) < t, 0 \leq i < n\},$$
$$F^*(f, x, y, t) = \limsup_{n\to\infty} \frac{1}{n} \#\{i \lceil d(f^i(x), f^i(y)) < t, 0 \leq i < n\}.$$
f is called distributional chaos if there exists an uncountable subset $S \subseteq X$ such that for every pair $x, y \in S$ of distinct points we have $F(f, x, y, \varepsilon) = 0 \ for \ some \ \varepsilon > 0 \ and \ F^*(f, x, y, t) = 1 \ for \ all \ t > 0$.

**Remark 1.** Distributional chaos was generalized in [6], so the above f is also called distributional chaos of type 1, briefly, DC1. If there exists an uncountable subset $S \subseteq X$ such that for every pair $x, y \in S$ of distinct points, we have
$$F(f, x, y, \varepsilon) < 1 \ for \ some \ \varepsilon > 0 \ and \ F^*(f, x, y, t) = 1 \ for \ all \ t > 0,$$
or
$$F(f, x, y, t) < F^*(f, x, y, t) \ for \ all \ t \in J, where \ J \ is \ some \ nondegenerate \ interval,$$
then we say that f exhibits distributional chaos of type 2-3, briefly, DC2 or DC3, respectively.

In the following definition, we will generalize further the distributional chaos. It is similar to DC2, so we call it DC2'.

**Definition 3.** Let $(X, d)$ be a metric space, and $f: X \to X$ be a continuous map. Note

$$F(f, x, y, t) = \liminf_{n \to \infty} \frac{1}{n} \#\{i \lceil d(f^i(x), f^i(y)) < t, 0 \leq i < n\},$$

$$F^*(f, x, y, t) = \limsup_{n \to \infty} \frac{1}{n} \#\{i \lceil d(f^i(x), f^i(y)) < t, 0 \leq i < n\},$$

f is called distributional chaos of type 2', briefly DC2', if there exists an uncountable subset $S \subseteq X$ such that for every pair $x, y \in S$ of distinct points we have $F(f, x, y, \varepsilon) = 0 \ for \ some \ \varepsilon > 0 \ and \ F^*(f, x, y, t) > 0 \ for \ all \ t > 0$.

**Definition 4.** Let $(X, d)$ be a metric space, $f: X \to X$ be a continuous map and $\{q_i\}$ be a strictly infinitely increasing sequence of positive integers. Note

$$F(f, x, y, t, q_i) = \liminf_{n \to \infty} \frac{1}{n} \#\{i \lceil d(f^{q_i}(x), f^{q_i}(y)) < t, 0 \leq i < n\},$$

$$F^*(f, x, y, t, q_i) = \limsup_{n \to \infty} \frac{1}{n} \#\{i \lceil d(f^{q_i}(x), f^{q_i}(y)) < t, 0 \leq i < n\},$$

f is called distributional chaos in a sequence briefly SDC if there exists an uncountable subset $S \subseteq X$ such that for every pair $x, y \in S$ of distinct points we have
$F(f, x, y, \varepsilon, q_i) = 0 \ for \ some \ \varepsilon > 0 \ and \ F^*(f, x, y, t, q_i) = 1 \ for \ all \ t > 0$.

**Remark 2.** We can also generalize the distributional chaos in a sequence named distributional chaos in a sequence of type1-3 which is similar to the distributional chaos, briefly SDC1-3. That is if
$F(f, x, y, \varepsilon, q_i) < 1 \ for \ some \ \varepsilon > 0 \ and \ F^*(f, x, y, t, q_i) = 1 \ for \ all \ t > 0$,
or
$F(f, x, y, t, q_i) < F^*(f, x, y, t, q_i) for \ all \ t \in J, where \ J \ is \ some \ nondegenerate \ interval$,
then we say that f exhibits distributional chaos in a sequence of type 2-3, briefly, SDC2 or SDC3, respectively.

From the above definitions it is easy to see that the following statements hold:
(a) DC1 implies DC2, SDC and DC2',
(b) both DC2 and DC2' respectively implies DC3,
(c) all DC1, DC2, DC2' and SDC respectively implies chaos in the sense of Li and Yorke,
(d) SDC1 implies SDC2, and SDC2 implies SDC3.

It has been proved that three versions of distributional chaos of type 1-3(DC1–DC3) are mutually nonequivalent[12] and DC1 is not equivalent to SDC1[11]. In the following example, we will show that DC2' is also neither equivalent to DC1 nor DC2. DC2' is not equivalent to DC3 by (c) and [13]. So the four versions of distributional chaos(DC1,DC2,DC2',DC3) are mutually nonequivalent.

**Example1**
Let $X = [0, +\infty)$ and define the metric d: $X \times X \to [0,1]$ be

$$d(x,y) = \begin{cases} 0, x = y \\ \frac{1}{2^k}, [x] = [y] \equiv 0 \pmod 2, x \neq y \in \left[\sum_{j=1}^{2k} b_j, \sum_{j=1}^{2k+1} b_j\right), k \in N, \\ 1, else \end{cases}$$

where $b_1 = 1, b_i = 2^{b_1+b_2+\cdots+b_{i-1}}$, it is easily to see that $(X, d)$ is a discrete metric space. Let $f: ([0, +\infty), d) \to ([0, +\infty), d)$ be $f(x) = x + 1$. We will claim that f is DC2'.

(i) Take $D_0 = (0,1), \forall 0 < t < \frac{1}{2}, \exists k_0 \in N$, such that $\frac{1}{2^{k_0}} < t$, note $\sum_{j=1}^{k} b_j = L_k. \forall x, y \in D_0 (x \neq y), \forall k \geq k_0, \forall 2i \in (L_{2k}, L_{2k+1} - 1)$, notices that $f^{2i}(x) = 2i + x \in (L_{2k}, L_{2k+1}), f^{2i}(y) = 2i + y \in (L_{2k}, L_{2k+1})$ and $[f^{2i}(x)] = [f^{2i}(y)] = 2i$, so $d\left(f^{2i}(x), f^{2i}(y)\right) = d(2i + x, 2i + y) = \frac{1}{2^k} \leq \frac{1}{2^{k_0}} < t$,

$$\frac{1}{L_{2k+1}} \#\{i \mid d(f^i(x), f^i(y)) < t, 0 \leq i < L_{2k+1}\} \geq \frac{1}{L_{2k+1}} \left(\frac{L_{2k+1} - 1 - L_{2k}}{2} - 1\right)$$

So that

$\limsup_{n \to \infty} \frac{1}{n} \#\{i \mid d(f^i(x), f^i(y)) < t, 0 \leq i < n\} \geq \limsup_{k \to \infty} \frac{1}{L_{2k+1}} \#\{i \mid d(f^i(x), f^i(y)) < t, 0 \leq i < L_{2k+1}\} \geq \limsup_{k \to \infty} \frac{1}{L_{2k+1}} \left(\frac{L_{2k+1}-1-L_{2k}}{2} - 1\right) = \limsup_{k \to \infty} \frac{\frac{2^{b_1+b_2+\cdots+b_{2k}}}{2} - 2}{b_1+b_2+\cdots+b_{2k}+2^{b_1+b_2+\cdots+b_{2k}}} = \frac{1}{2}$.

(1)

Since $[x] = [y] \equiv 0 \pmod 2$ where $x, y \in D_0 (x \neq y)$, then $\forall i \in N, [f^{2i+1}(x)] = [x] + 2i + 1 = [y] + 2i + 1 = [f^{2i+1}(y)] \equiv 1 \pmod 2$, so $d\left(f^{2i+1}(x), f^{2i+1}(y)\right) = 1$.

So that $\limsup_{n \to \infty} \frac{1}{n} \#\{i \mid d(f^i(x), f^i(y)) < t, 0 \leq i < n\} \leq \limsup_{n \to \infty} \frac{\left[\frac{n}{2}\right]+1}{n} = \frac{1}{2}$. (2)

Therefore, $\limsup_{n \to \infty} \frac{1}{n} \#\{i \mid d(f^i(x), f^i(y)) < t, 0 \leq i < n\} = \frac{1}{2}$ by (1) and (2).

(ii) $\forall x, y \in D_0 (x \neq y), \forall k \in N, \forall i \in (L_{2k+1}, L_{2k+2} - 1)$, notices that

$$f^i(x) \neq f^i(y) \in (L_{2k+1}, L_{2k+2}), so\ d\left(f^i(x), f^i(y)\right) = 1.$$

$$\frac{1}{L_{2k+2}} \#\{i \mid d(f^i(x), f^i(y)) < \frac{1}{2}, 0 \leq i < L_{2k+2}\} \leq \frac{1}{L_{2k+2}} (L_{2k+1} + 1)$$

$\liminf_{n \to \infty} \frac{1}{n} \#\{i \mid d(f^i(x), f^i(y)) < \frac{1}{2}, 0 \leq i < n\} \leq \liminf_{k \to \infty} \frac{1}{L_{2k+2}} \#\{i \mid d(f^i(x), f^i(y)) < \frac{1}{2}, 0 \leq i < L_{2k+2}\} \leq \liminf_{k \to \infty} \frac{1}{L_{2k+2}} (L_{2k+1} + 1) = \liminf_{k \to \infty} \frac{b_1+b_2+\cdots+b_{2k+1}+1}{b_1+b_2+\cdots+b_{2k+1}+2^{b_1+b_2+\cdots+b_{2k+1}}} = 0$.

Therefore f is DC2' by (i) and (ii).

**Definition 5.** Let $\omega(x, f)$ denote the set of $\omega$-limit points of f, and $S \subseteq X$. We say that S is an $\omega$-scrambled set if, for any $x, y \in S$ with $x \neq y$,

(1) $\omega(x, f) \setminus \omega(x, f)$ is uncountable;
(2) $\omega(x, f) \cap \omega(x, f)$ is nonempty; and
(3) $\omega(x, f)$ is not contained in the set of periodic points.

We say that f is $\omega$-chaotic, if there exists an uncountable $\omega$-scrambled set.

**Definition 6.** We say that a semiflow (T, X) is Ruelle–Takens chaotic (R-T chaos) if it is point-transitive and sensitive.

**Lemma 1**[7] Let f: I → I be a continuous map, then f is Li-Yorke chaos if and only if f is distributional chaos in a sequence.

**Lemma 2**[4] Let f: I → I be a continuous map, then f is ω − chaos if and only f has positive topological entropy.

**Lemma 3**[14] Let X be a compact metric space, f: X → X be continuous, x, y ∈ X, ∀N > 0, we have
(i) If for t > 0, $F(f, x, y, t) = 0$, then $F(f^N, x, y, t) = 0$.
(ii) If for t > 0, $F(f^N, x, y, t) = 0$, then $F(f, x, y, t) = 0$.

## 3. Main results

**Theorem 1.** Let $(X, f)$ *be compact* dynamical system and $Q = \{q_i\}$ be a strictly infinitely increasing sequence of positive integers. If ∃M such that ∀i, we have $q_{i+1} - q_i \leq M$, then that f is distributional chaos if and only if that f is distributional chaos in a sequence Q.

**Proof.** Necessity. (i) Since f is distributional chaos, $F(f, x, y, \varepsilon) = 0$, then there is an increasing sequence $\{n_k\}$ of positive integers such that

$$\text{for } k \to \infty, \frac{1}{n_k} \#\{i \mid d\left(f^i(x), f^i(y)\right) < \varepsilon, 0 \leq i < n_k\} \to 0. \tag{3}$$

Put $m_k = \left[\frac{n_k}{M}\right]$.

Then for each k,

$$\#\{i \mid d(f^{q_i}(x), f^{q_i}(y)) < \varepsilon, 0 \leq i < m_k\} \leq \#\{i \mid d\left(f^i(x), f^i(y)\right) < \varepsilon, 0 \leq i < n_k\},$$

and further

$$\frac{1}{n_k}\#\{i \mid d(f^{q_i}(x), f^{q_i}(y)) < \varepsilon, 0 \leq i < m_k\} \leq \frac{1}{n_k}\#\{i \mid d\left(f^i(x), f^i(y)\right) < \varepsilon, 0 \leq i < n_k\}.$$

This gives for $k \to \infty, \frac{F(f, x, y, \varepsilon, q_i)}{M} = 0$ by (3). Therefore $F(f, x, y, \varepsilon, q_i) = 0$.

(ii) Since f is distributional chaos, $F^*(f, x, y, t) = 1$, for all t > 0, then there is an increasing sequence $\{n_k\}$ of positive integers such that

$$\text{for } k \to \infty, \frac{1}{n_k}\#\{i \mid d\left(f^i(x), f^i(y)\right) < t, 0 \leq i < n_k\} \to 1,$$

and so

$$\frac{1}{n_k}\#\{i \mid d\left(f^i(x), f^i(y)\right) \geq t, 0 \leq i < n_k\} \to 0. \tag{4}$$

Put $m_k = \left[\frac{n_k}{M}\right]$.

Then for each k,

$$\#\{i \mid d(f^{q_i}(x), f^{q_i}(y)) \geq t, 0 \leq i < m_k\} \leq \#\{i \mid d(f^i(x), f^i(y)) \geq t, 0 \leq i < n_k\}$$

And further

$$\frac{1}{n_k}\#\{i \mid d(f^{q_i}(x), f^{q_i}(y)) \geq t, 0 \leq i < m_k\} \leq \frac{1}{n_k}\#\{i \mid d(f^i(x), f^i(y)) \geq t, 0 \leq i < n_k\}.$$

This gives for $k \to \infty$, $\frac{1}{M \cdot m_k}\#\{i \mid d(f^{q_i}(x), f^{q_i}(y)) \geq t, 0 \leq i < m_k\} \to 0$ by (4).

Therefore $\frac{1}{m_k}\#\{i \mid d(f^{q_i}(x), f^{q_i}(y)) \geq t, 0 \leq i < m_k\} \to 0$ for $k \to \infty$,

$$\frac{1}{m_k}\#\{i \mid d(f^{q_i}(x), f^{q_i}(y)) < t, 0 \leq i < m_k\} \to 1 \text{ for } k \to \infty,$$

$F^*(f, x, y, t, q_i) = 1$.

Therefore f is distributional chaos in a sequence Q.

Sufficiency.(iii) Since f is distributional chaos in a sequence Q, $F(f, x, y, \varepsilon, q_i) = 0$, then there is an increasing sequence $\{n_k\}$ of positive integers such that for $k \to \infty$, $\frac{1}{n_k}\#\{i \mid d(f^{q_i}(x), f^{q_i}(y)) < \varepsilon, 0 \leq i < n_k\} \to 0$ (5)

Since $f$ is continuous and X is compact, $f^j$ is uniformly continuous for each $j = 1, 2, \ldots, M$. Consequently, for fixed $\varepsilon > 0$, there exists $s > 0$ such that

$d(f^{q_i}(x), f^{q_i}(y)) < \varepsilon$ when $d(f^{q_i-j+1}(x), f^{q_i-j+1}(y)) < s$, for each $j = 1, 2, \ldots, M$,

so

$$\sum_{j=1}^{M} \#\{i \mid d(f^{q_i-j+1}(x), f^{q_i-j+1}(y)) < s, 0 \leq i < n_k\} \leq M \cdot \#\{i \mid d(f^{q_i}(x), f^{q_i}(y)) < \varepsilon, 0 \leq i < n_k\}.$$

Since $q_{i+1} - q_i \leq M$,

$$\frac{1}{n_k}\#\{i \mid d(f^i(x), f^i(y)) < s, 0 \leq i < n_k\} \leq \frac{1}{n_k}\sum_{j=1}^{M}\#\{i \mid d(f^{q_i-j+1}(x), f^{q_i-j+1}(y)) < s, 0 \leq i < n_k\} \leq M \cdot \frac{1}{n_k}\#\{i \mid d(f^{q_i}(x), f^{q_i}(y)) < \varepsilon, 0 \leq i < n_k\}$$ (6)

This shows that $F(f, x, y, s) = 0$ for $k \to \infty$ by (5) and (6).

(iv) Since f is distributional chaos in a sequence Q, $F^*(f, x, y, p, q_i) = 1$ for all $p > 0$, then there is an increasing sequence $\{n_k\}$ of positive integers such that for $k \to \infty$, $\frac{1}{n_k}\#\{i \mid d(f^{q_i}(x), f^{q_i}(y)) < p, 0 \leq i < n_k\} \to 1$,

and so

$$\frac{1}{n_k}\#\{i \mid d(f^i(x), f^i(y)) \geq p, 0 \leq i < n_k\} \to 0.$$ (7)

Since $f$ is continuous and X is compact, $f^j$ is uniformly continuous for each $j = 1, 2, \ldots, M$. Consequently, for $\forall t > 0$, there exists $s > 0$ such that

$d(f^{q_i+j-1}(x), f^{q_i+j-1}(y)) < t$ when $d(f^{q_i}(x), f^{q_i}(y)) < s$, for each $j = 1, 2, \ldots, M$.

So we have $d(f^{q_i}(x), f^{q_i}(y)) \geq s$ when $d(f^{q_i+j-1}(x), f^{q_i+j-1}(y)) \geq t$ for each $j = 1, 2, \ldots, M$.

So

$$\sum_{j=1}^{M} \#\{i \mid d\left(f^{q_i+j-1}(x), f^{q_i+j-1}(y)\right) \geq t, 0 \leq i < n_k\} \leq M \cdot \#\{i \mid d(f^{q_i}(x), f^{q_i}(y)) \geq s, 0 \leq i < n_k\}$$

Since $q_{i+1} - q_i \leq M$,

$$\frac{1}{n_k}\#\{i \mid d\left(f^i(x), f^i(y)\right) \geq t, 0 \leq i < n_k\} \leq \frac{1}{n_k}\sum_{j=1}^{M} \#\{i \mid d\left(f^{q_i+j-1}(x), f^{q_i+j-1}(y)\right) \geq t, 0 \leq i < n_k\} + \frac{q_0}{n_k} \leq M \cdot \frac{1}{n_k}\#\{i \mid d(f^{q_i}(x), f^{q_i}(y)) \geq s, 0 \leq i < n_k\} + \frac{q_0}{n_k}. \quad (8)$$

This shows that $\frac{1}{n_k}\#\left\{i \mid d\left(f^i(x), f^i(y)\right) \geq t, 0 \leq i < n_k\right\} \to 0$ for $k \to \infty$ by (7) and (8).

Therefore

$$\frac{1}{n_k}\#\{i \mid d\left(f^i(x), f^{q_i}(y)\right) < t, 0 \leq i < n_k\} \to 1 \text{ for } k \to \infty$$

This means that $F^*(f, x, y, t) = 1$.

Therefore f is distributional chaos.

**Remark 3.** The theorem1 is the generalization of [14], where N=M, i.e. $Q = \{q_i\} = \{Ni\}$. If f is distributional chaos in a sequence $Q = \{q_i\}$ but is not distributional chaos, then we have $\limsup_{i \to \infty} q_{i+1} - q_i = \infty$.

**Theorem 2.** Let X be a compact metric space. Then f is DC2' if and only if $f^N$ i is too.

**Proof.** Necessity. Sincef is DC2' ,so $F^*(f, x, y, t) > 0$ for all $t > 0$, then there exists an increasing sequence $\{m_k\}$ of positive integers such that

$$\text{for all } t > 0, \frac{1}{m_k}\#\{i \mid d\left(f^i(x), f^i(y)\right) < t, 0 \leq i < m_k\} > 0 \ (k \to \infty) \quad (9)$$

Since X is compact, $f^i$ is uniformly continuous for each $i = 1, 2, \ldots, N$. Consequently, $\forall s > 0$, there exists $p > 0$ such that for all u,v and each $i = 1, 2, \ldots, N, d\left(f^i(u), f^i(v)\right) \geq p$ whenever $d(f^N(u), f^N(v)) \geq s$. This imply

$$N(\#\{i \mid d\left(f^{Ni}(x), f^{Ni}(y)\right) \geq s, 0 \leq i < n_k\} - 1)$$

$$\leq \#\{i \mid d\left(f^i(x), f^i(y)\right) \geq p, 0 \leq i < N \cdot n_k\}$$

(10).

Take $n_k = [\frac{m_k}{N}]$ for all $k \geq 1$, by a simple calculation, we may derive from (10) that

$$\frac{1}{m_k}\#\{i \mid d\left(f^i(x), f^i(y)\right) < p, 0 \leq i < m_k\} \leq \frac{1}{n_k}\#\{i \mid d\left(f^{Ni}(x), f^{Ni}(y)\right) < s, 0 \leq i < n_k\} + \frac{1}{n_k}. \quad (11)$$

By (9)and(11), we have

$$\text{for all } s > 0, \frac{1}{n_k}\#\{i \mid d\left(f^{Ni}(x), f^{Ni}(y)\right) < s, 0 \leq i < n_k\} > 0 \ (k \to \infty). \quad (12)$$

(12)and lemma 3 show that $f^N$ is DC2'.

Sufficiency. Since $f^N$ is DC2', so $F^*(f^N, x, y, t) > 0$ $for\ all\ t > 0$, then there exists an increasing sequence $\{m_k\}$ of positive integers such that

$for\ all\ t > 0$ , $\frac{1}{m_k}\#\{i \mid d\left(f^{Ni}(x), f^{Ni}(y)\right) < t, 0 \leq i < m_k\} > 0$ ( k → ∞)

Take $n_k = m_k \cdot N$ ,
then for all k ≥ 1 and $all\ t > 0$,

$$\frac{1}{m_k}\#\{i \mid d\left(f^i(x), f^i(y)\right) < t, 0 \leq i < n_k\}$$
$$\geq \frac{1}{m_k}\#\{i \mid d\left(f^{Ni}(x), f^{Ni}(y)\right) < t, 0 \leq i < m_k\} > 0,\ (\text{k} \to \infty)$$

so $\frac{1}{n_k}\#\{i \mid d\left(f^i(x), f^i(y)\right) < t, 0 \leq i < n_k\} > 0(\text{k} \to \infty)$. (13)

(13) and lemma 3 show that f is DC2'.

**Theorem 3.   SDC and ω−chaos don't imply each other.**
**Proof.** On one hand, let f: I → I be a continuous map, and be Li-Yorke chaos with zero topological entropy , then f is distributional chaos in a sequence by lemma 1, but is not ω −chaos by lemma 2.

On the other hand, there is a dynamical system $(X, g)$ that g is ω −chaos but is not Li-Yorke chaos in [15]. Because distributional chaos in a sequence implies Li-Yorke chaos, g is not distributional chaos in a sequence.

**Theorem 4.   SDC and R-T chaos don't imply each other.**

**Proof.** On one hand, let $(X, f)$ be a minimal distal non equicontinuity dynamical system, then $(X, f)$ is R-T chaos but not Li-Yorke chaos, so $(X, f)$ is not distributional chaos in a sequence like the proof of theorem 3.

On the other hand , Professor Wang and Li give a dynamical system $(X, f)$ , which is distributional chaos in a sequence but not R-T chaos in[16].

**Theorem 5.   SDC and DC3 don't imply each other.**
**Proof.** On one hand, there is a continuous map f: I → I, which is Li-Yorke chaos but not DC3 in[3], then f is distributional chaos in a sequence by lemma1.
On the other hand, there is a continuous map f which is DC3 but not Li-Yorke chaos in [6], then f is not distributional chaos in a sequence like the proof of theorem3.